\documentclass[11pt,reqno]{amsart}
\usepackage{amsfonts}
\usepackage{fancyhdr}
\usepackage{times}
\usepackage[T1]{fontenc}
\usepackage{mathrsfs}
\usepackage{latexsym}
\usepackage[dvips]{graphics}
\usepackage{epsfig}
\usepackage{hyperref, amsmath, amsthm, amsfonts, amscd, flafter,epsf}
\usepackage{amsmath,amsfonts,amsthm,amssymb,amscd}
\usepackage{color}
\usepackage[parfill]{parskip}
 \usepackage{mathrsfs} 
\usepackage{geometry}                
\usepackage{epstopdf}

\usepackage{verbatim}

\newtheorem{thm}{Theorem}

\newtheorem{lem}[thm]{Lemma}
\newtheorem{thm*}{Theorem}
\newtheorem*{con*}{Conjecture}
\newtheorem*{lem*}{Lemma}

\def\g{\gamma}

\begin{document}
 \baselineskip=17pt
\hbox{}
\medskip
\title{The Number of Zeros of $\zeta'(s)$}
\author{Fan Ge}

\email{fan.ge@rochester.edu}

\address{Department of Mathematics, University of Rochester, Rochester, NY}

\thanks{This work was partially supported by NSF grant DMS-1200582.}

 \maketitle

\begin{abstract}
Assuming the Riemann Hypothesis, we prove that
$$
N_1(T) = \frac{T}{2\pi}\log \frac{T}{4\pi e} + O\bigg(\frac{\log
T}{\log\log T}\bigg),
$$
where $N_1(T)$ is the number of zeros of $\zeta'(s)$ in the region
$0<\Im s\le T$.
\end{abstract}

\section{Introduction}

The distribution of zeros of the first derivative of the Riemann
zeta-function is interesting and important  both in its own right and because of its connection with
the distribution of zeros of zeta. As just one illustration of the latter we cite
A. Speiser's~\cite{Spe} theorem that the Riemann Hypothesis is
equivalent to the nonexistence of non-real zeros of $\zeta'(s)$ in
the half-plane $\Re s<1/2$.

Let $\rho'=\beta'+i\g'$ denote a  generic nonreal zero of  $\zeta'(s)$,
where $s=\sigma+it$ is a complex variable,
and for $T\geq2$ consider the zero counting function
$$
N_1(T):=\sum_{\substack{0<\g'\leq T\\ \beta'>0}} 1.
$$
In~\cite{Ber}  B. C. Berndt proved that
\begin{align}
\label{eq N1}
N_1(T)=\frac{T}{2\pi}\log \frac{T}{4\pi e} + O(\log T)
\end{align}
as $T\to\infty$.
This should   be compared with the well-known formula (see  \cite{Tit})
\begin{align}
\label{eq N}
N(T)=\frac{T}{2\pi}\log \frac{T}{2\pi e} + O(\log T),
\end{align}
which counts the zeros $\rho=\beta+i\g$ of $\zeta(s)$ in the same region. That is,
\begin{equation}\label{ }\notag
N(T)=\sum_{\substack{0<\g\leq T\\ \beta>0}} 1.
\end{equation}
By an old result of Littlewood~\cite{Lit}, 
if the Riemann Hypothesis (RH) is true,   the error term in \eqref{eq N} may be replaced by \begin{equation}\label{Lit Er}
O\bigg( \frac{\log  T}{\log\log T}\bigg).
\end{equation}
 Recently  H. Akatsuka~\cite{Aka} showed that  if RH is true,   the error
term in \eqref{eq N1}   may  also be reduced, namely  to
\begin{equation}\label{Aka Er}
O\bigg(\frac{\log T}{\sqrt{\log\log T}}\bigg).
\end{equation}

Our purpose here is to  show that if RH is true, the   error term \eqref{Lit Er} holds
in both  \eqref{eq N1}  and \eqref{eq N}.

\begin{thm}\label{thm N1}
Assume RH. Then we have
$$
N_1(T) = \frac{T}{2\pi}\log \frac{T}{4\pi e} + O\bigg(\frac{\log
T}{\log\log T}\bigg)
$$
as $T\to\infty$.
\end{thm}

D. W. Farmer, S. M. Gonek and C. P. Hughes \cite{FGH} have
conjectured that the error term in \eqref{eq N} is $O(\sqrt{\log T
\log\log T})$. This raises the question of what error term one
should expect in  \eqref{eq N1} for a given error term in \eqref{eq
N}. By slightly modifying our proof of Theorem \ref{thm N1}, we can
show
\begin{thm}\label{thm general}
Assume RH and suppose that the error term in \eqref{eq N} is
$O(\Phi(T))$ for some increasing function $\log\log T \ll \Phi(T)\ll
\log T$. Then we have
\begin{equation}\label{N1 thm}
 N_1(T) = \frac{T}{2\pi}\log \frac{T}{4\pi e} +
O\big(\max\big\{\Phi(2T), \sqrt{\log T}\log\log T \big\}\big).
\end{equation}
\end{thm}
Using this with the conjecture of Farmer et al., we obtain
 $$ N_1(T) = \frac{T}{2\pi}\log \frac{T}{4\pi e} +
O\big(  \sqrt{\log T}\log\log T \big).
$$
Of course, we really expect the error term in \eqref{N1 thm} to be
$O(\Phi(T))$ in general.

\section{lemmas}

We first state three results of Akatsuka.

\begin{lem}\label{lem N1=arg..}
Assume RH. For $T\ge 2$ satisfying $\zeta(\sigma+iT)\ne 0$ and
$G(\sigma+iT)\ne 0$ for any $\sigma\in\mathbb{R}$, we have
$$
N_1(T) = \frac{T}{2\pi}\log \frac{T}{4\pi e} + \frac{1}{2\pi}\arg G(1/2+iT)+\frac{1}{2\pi}\arg \zeta(1/2+iT) +O(1),
$$
where $$G(s)=\frac{-2^s}{\log 2}\zeta'(s),$$ and the argument is
defined by continuous variation from $+\infty$, with the argument at
$+\infty$ being 0.
\end{lem}
This is Proposition 3.1 in \cite{Aka}.

\begin{lem}\label{lem arg far away}
Assume RH. For $1/2+(\log\log T)^2/\log T \le \sigma \le 3/4$, we
have
$$\arg G(\sigma +iT) \ll \frac{(\log T)^{2(1-\sigma)}}{\log\log T}.$$
\end{lem}
  See Remark 2.5 in \cite{Aka}.

\begin{lem}\label{lem arg g over zeta}
Assume RH. Then for $1/2<\sigma<20$, we have
$$
\arg\bigg(-\frac{2^{\sigma+iT}}{\log
2}\frac{\zeta'}{\zeta}(\sigma+iT)\bigg) =O\bigg(\frac{\log\log
T}{\sigma-\frac{1}{2}}\bigg)\ .
$$
\end{lem}
See Lemma 2.3 in \cite{Aka}.

We also require the following lemma.
\begin{lem}\label{lem zeta'}
For all $t$ sufficiently large we have
$$
\frac{\zeta''}{\zeta'}(s)=\sum_{|\rho'-s|<5}\frac{1}{s-\rho'}
+O(\log
t),
$$
uniformly for $-1\le \sigma \le 2$.
\end{lem}
 This can be proved in a standard way. See Theorem 9.6 (A) in
\cite{Tit} for example.

Let us define
$$
F_1(t)=\sum_{\beta'>1/2}
\frac{\beta'-1/2}{(\beta'-1/2)^2+(\gamma'-t)^2}.
$$
We require the following result of Y. Zhang \cite{Zha}.

\begin{lem}\label{lem F1}
Assume RH and order  the ordinates of the zeros of  $\zeta(s)$ as $0<\g_1\leq \g_2\leq \cdots$.
Then
$$
\int_{\g_n}^{\g_{n+1}} F_1(t) \ll 1,
$$
where the implied constant is absolute.
\end{lem}
This is a combination of Lemma 4 in \cite{Zha} and equation
(4.1) in \cite{Zha}.

To state our final lemma we
write
$$
H=\frac{(\log\log T)^3}{\log T}
$$
and let $N_d$ denote the number of distinct zeros of $\zeta(s)$ on the
vertical segment
$$
\big [1/2 + i (T-H), 1/2+i (T+H)\big].
$$
We also let  $N_{1}(\mathscr R)$ be the number of zeros of
$\zeta'(s)$ in the rectangular region $\mathscr{R}$ given by
 $$
 T-H\le t \le T+H,\ \ \ \
 1/2<\sigma \le 1/2+H.
 $$

\begin{lem}\label{lem Nd}
Assume RH. We have $$N_{1}(\mathscr R)\ll N_d + 1.$$
\end{lem}
\proof  Consider the integral
$$
\mathcal{I}:=\int_{T-H}^{T+H} F_1(t)dt.
$$
From Lemma \ref{lem F1} we see that
\begin{equation}\label{N_d}
\mathcal{I}\ll N_d + 1.
\end{equation}
Thus, to prove the lemma, it suffices to show that
$$N_{1}(\mathscr R)\ll \mathcal{I}.$$
First observe that
$$
F_1(t)\ge \sum_{\rho'\in
 \mathscr R}\frac{\beta'-1/2}{(\beta'-1/2)^2+(\gamma'-t)^2} \ ,
$$
since each summand is   positive. Hence,
\begin{align*}
\mathcal{I} & \ge \int_{T-H}^{T+H}
\sum_{\rho'\in \mathscr{R}}\frac{\beta'-1/2}{(\beta'-1/2)^2+(\gamma'-t)^2} dt
\\
& = \sum_{\rho'\in \mathscr{R}} \int_{T-H}^{T+H}
\frac{\beta'-1/2}{(\beta'-1/2)^2+(\gamma'-t)^2} dt.
\end{align*}
Let $ \theta(\rho')\in (0, \pi)$ be the argument of the angle at
$\rho'$ with two rays through $1/2+i(T-H)$ and $1/2+i(T+H)$
respectively. It is easy to see that
$$
\theta(\rho')=\int_{T-H}^{T+H}
\frac{\beta'-1/2}{(\beta'-1/2)^2+(\gamma'-t)^2} dt.
$$
This gives us
$$
 \mathcal{I}  \ge  \sum_{\rho'\in \mathscr{R}} \theta(\rho').
$$
Now for $\rho'\in \mathscr{R}$, we clearly have
$\theta(\rho')\ge c$ for some absolute positive constants $c$.
Therefore, we have
$$\mathcal{I}\ge c N_{1}(\mathscr R)\gg N_{1}(\mathscr R).$$
The result follows on combining this  and \eqref{N_d}.
  \qed

\section{proof of Theorem \ref{thm N1}}

Since $N_1(T)$ is right
continuous with respect to $T$, it suffices to consider $T$'s such that
$\zeta(\sigma+iT)\ne 0$ and $G(\sigma+iT)\ne 0$ for any
$\sigma\in\mathbb{R}$.
It is well known (see \cite{Tit}) that RH implies
$$
\arg\zeta(1/2+iT)\ll \frac{\log T}{\log\log T},
$$
(in fact, this is equivalent to  \eqref{Lit Er}).
Thus, in view of Lemma \ref{lem N1=arg..}, to prove the theorem it
suffices to prove that
\begin{align}\label{bound arg g}
 \arg G(1/2+iT)\ll \frac{\log T}{\log\log T}.
\end{align}
Recall that $ \arg G(1/2+iT)$ is defined by continuous variation along the horizontal
line from $\infty +iT$ to $1/2+iT$ starting with the value $0$.
Let $ \Delta_1$ be the change in  argument of $G$ along the horizontal line
from $\infty+iT$ to $1/2 + (\log\log T)^2/\log T +iT$, and let
$ \Delta_2$  be  the change  along the horizontal segment from $1/2 + (\log\log
T)^2/\log T +iT$ to $1/2+iT$.  We shall show that  $ \Delta_1$ and $ \Delta_2$ are
both bounded
by $O(\log T/\log\log T)$.

The estimation of $ \Delta_1$ is immediate, for by Lemma \ref{lem arg far away}
with $\sigma= 1/2+ (\log\log T)^2/\log T $ we see
that
$$
\Delta_1 = \arg G(\sigma +iT) \ll \frac{(\log T)^{2(1-\sigma )}}{\log\log
T}\ll \frac{\log T}{\log\log T}.
$$

Next we bound  $ \Delta_2$.  Using the notation
$$
H=\frac{\log\log^3 T}{\log T},
$$
 we see that
$$
\Delta_2= \Im \int_{1/2}^{1/2+ H/\log\log T} \frac{G'}{G}(\sigma+iT)d\sigma.
$$
By the definition of $G(s)$ we have $(G'/G)(s)=(\zeta''/\zeta')(s)+O(1)$. Thus,
it follows from Lemma \ref{lem zeta'}  that
$$
\frac{G'}{G}(s)=\sum_{|\rho'-s|<5}\frac{1}{s-\rho'}+O(\log T)
$$
for  $s$  on the  segment $[1/2 +iT, 1/2+ H/\log\log  T +iT]$.
It  is convenient to modify this formula slightly.
Let $\mathscr{D}$ be the  disk centered at $1/2+
 {H}/{2\log\log  T} +iT$ with radius $5$.
We  claim that
 $$
\frac{G'}{G}(s)=\sum_{\rho'\in \mathscr{D}}\frac{1}{s-\rho'}+O(\log
T)
$$
for $s$  on the horizontal segment $[1/2 +iT, 1/2+ H/\log\log
T +iT]$.
 Indeed, for such $s$ and for $\rho'$ not belonging to the
 intersection of $\mathscr{D}$ and the disk $|\rho'-s|<5$,
 we see that $(s-\rho')^{-1}\ll1$.
 Since the number of such $\rho'$ is $O(\log
 T)$, their contribution   is also  $O(\log T)$.


It now follows that
\begin{align}\label{Delta_2}
\Delta_2 & = \Im \int_{1/2+iT}^{1/2+ H/\log\log
T+iT}\bigg(\sum_{\rho'\in \mathscr{D}}\frac{1}{s-\rho'}+O(\log T)\bigg)ds  \notag\\
& =\sum_{\rho'\in  \mathscr{D}} \bigg(\Im\int_{1/2+iT}^{1/2+ H/\log\log
T+iT}\frac{1}{s-\rho'}\;ds\bigg)+O((\log\log T)^2)\\
& =\sum_{\rho'\in  \mathscr{D}} \big[ \arg (1/2 + H/\log\log  T
+iT - \rho')  -  \arg (1/2
+iT - \rho')  \big] +O((\log\log T)^2) \notag \\
& = \sum_{\rho'\in  \mathscr{D}} f(\rho')\ +\ O((\log\log T)^2) ,\notag
\end{align}
say. Notice that $ f(\rho')$ is plus or minus the argument of the angle subtended by the
segments from $1/2+iT$ to $\rho'$ and from $1/2+H/\log\log T+iT$ to $\rho'.$
Thus, in particular, $ f(\rho')\ll 1$.

It remains to prove that
$$\sum_{\rho'\in \mathscr{D}} f(\rho') \ll \frac{\log T}{\log\log T}.$$
We split the sum into three parts. We let $\sum_1$ denote the  sum  over the $\rho'\in \mathscr{R}$, that is, the
$\rho'$ satisfying
$$
T-H\le \g'\le T+H, \ \ \ \ \ 1/2< \beta'\le 1/2+H.
$$
We let $\sum_2$ be the sum  over the zeros $\rho'=1/2+i\g'$, if any, with
$$
T-H\le \g' \le T+H.
 $$
Finally, we let $\sum_3$ denote the sum over  the remaining $\rho'$
in $\mathscr{D}$.

By Lemma \ref{lem Nd} and our observation above that  $f(\rho')\ll 1$, we
see that
$$
\sideset{}{_1}\sum \ll \max |f(\rho')|\cdot N_{1}(\mathscr R)
\ll  N_{1}(\mathscr R)\ll N_d
+ 1.
$$
Recall that RH implies that
$$
N(T)=\frac{T}{2\pi}\log \frac{T}{2\pi e} + O\Big(\frac{\log T}{\log\log T}\Big).
$$
Thus   $N_d\ll N(T+H)-N(T-H) \ll\log T/\log\log T$.
We therefore obtain
\begin{align}\label{eq sum1}
\sideset{}{_1}\sum \ll \frac{\log T}{\log\log T}.
\end{align}

A similar argument (together with the well-known fact that
$\zeta'(1/2+it)=0\Rightarrow \zeta(1/2+it)=0$) shows that
\begin{align}\label{eq sum2}
\sideset{}{_2}\sum \ll \frac{\log T}{\log\log T}.
\end{align}

Now we consider $\sum_3$. Recall that RH implies that $\zeta'(s)$
has no non-real zero to the left of 1/2. Thus, it is easy to see
that
$$
f(\rho')=\arg (1/2 + H/\log\log T +iT - \rho')  -  \arg (1/2 +iT -
\rho')\ll \frac{1}{\log\log T}
$$
for $\rho'$ in $\sum_3$. Since the number of $\rho'$ in $\mathscr D$ is
at most $O(\log T)$, we see that
$$
\sideset{}{_3}\sum \ll \frac{\log T}{\log\log T}.
$$

Combining this with (\ref{eq sum1}) and (\ref{eq sum2}), we see that
$$\sum_{\rho'\in \mathscr{D}} f(\rho') \ll \frac{\log T}{\log\log T}.
$$
Thus, by \eqref{Delta_2}   we obtain
 $$
 \Delta_2 \ll \frac{\log T}{\log\log T}.
 $$
This completes our proof. \qed


\section{Proof of Theorem 2}

Let $X$, $g$ and $k$ be positive parameters with $X=o(1)$, $g\ge 2$
and $k\in \mathbb{Z}$, to be determined later. As before, we let $
\Delta_1(X)$ be the change in argument of $G$ along the horizontal
line from $\infty+iT$ to $1/2 + X +iT$, and let $ \Delta_2(X)$ be
the change along the horizontal segment from $1/2 + X +iT$ to
$1/2+iT$.


It is standard (see Theorem 14.14 (B) in \cite{Tit}) to show that
$$
\arg \zeta (\sigma + iT)\ll \Phi(2T)
$$
for $\sigma\ge 1/2$. Using this with Lemma \ref{lem arg g over
zeta}, we obtain
\begin{align}\label{eq Delta 1 X}
\Delta_1(X)\ll \frac{\log\log T}{X} + \Phi(2T).
\end{align}

Next we bound  $ \Delta_2$.  For $0\le j\le k$ let
$$
Y_j=Xg^j,
$$
and define $\mathscr R_j$ to be the rectangular region
$$
\frac{1}{2}<\sigma\le \frac{1}{2}+Y_j,\ \ \ \ \ \  T-Y_j\le t \le
T+Y_j.
$$
Let $\mathscr U_1=\mathscr R_1$, and for $2\le j\le k$ let $\mathscr
U_j=\mathscr R_j-\mathscr R_{j-1}$. We denote by $N(\mathscr U_j)$
the number of zeros of $\zeta'(s)$ in $\mathscr U_j$.

Let $\mathscr{D}=\mathscr{D}(X)$ be the  disk centered at $1/2+
 X/2 +iT$ with radius $5$.
As before we have
 $$
\frac{G'}{G}(s)=\sum_{\rho'\in \mathscr{D}}\frac{1}{s-\rho'}+O(\log
T)
$$
for $s$  on the horizontal segment $[1/2 +iT, 1/2+X +iT]$. It
follows that
\begin{align}\label{Delta_2 X}
\Delta_2(X) = \sum_{\rho'\in  \mathscr{D}} f_X(\rho')\ +\ O(X\log T)
\end{align}
where $ f_X(\rho')$ is plus or minus the argument of the angle
subtended by the segments from $1/2+iT$ to $\rho'$ and from
$1/2+X+iT$ to $\rho'.$

We split the sum into $k+2$ parts. For $1\le j\le k$ we let $\sum_j$
denote the sum over the $\rho'\in \mathscr U_j$. We let $\sum_{k+1}$
be the sum  over the zeros $\rho'=1/2+i\g'$, if any, with
$$
T-Y_k\le \g' \le T+Y_k.
 $$
Finally, we let $\sum_{k+2}$ denote the sum over  the remaining
$\rho'$ in $\mathscr{D}$.

Using the same argument as in Theorem 1, we can show that
$$
\sideset{}{_j}\sum\ll g^{1-j} N(\mathscr U_j)
$$
for $1\le j\le k$. It follows from Lemma \ref{lem Nd} that
$$N(\mathscr U_j)\ll N(T+Y_j)-N(T-Y_j)+1\ll Y_j\log T + \Phi(T+Y_j).$$
This gives us
\begin{align*}
\sideset{}{_j}\sum \ll g^{1-j}(Y_j\log T + \Phi(T+Y_j))\le Xg\log T
+ g^{1-j}\Phi(T+Y_k)
\end{align*}
for $1\le j\le k$. Thus, we have
\begin{align*}
\sideset{}{_1}\sum+\sideset{}{_2}\sum+\cdots+\sideset{}{_k}\sum\ll
Xgk\log T + \Phi(T+Y_k).
\end{align*}
Similarly, we get
\begin{align*}
\sideset{}{_{k+1}}\sum\ll Xgk\log T + \Phi(T+Y_k),
\end{align*}
and
\begin{align*}
\sideset{}{_{k+2}}\sum\ll g^{-k}\log T.
\end{align*}
Combining the above estimates, we see that
$$\sum_{\rho'\in \mathscr{D}} f_X(\rho') \ll Xgk\log T + \Phi(T+Y_k)+g^{-k}\log
T.
$$
Thus, by \eqref{Delta_2 X}  we obtain
 $$
 \Delta_2 \ll Xgk\log T + \Phi(T+Y_k)+g^{-k}\log
T .
 $$
This together with \eqref{eq Delta 1 X} gives us
\begin{align}\label{eq bound Delta 2 X}
\arg G(1/2+iT)\ll  \frac{\log\log T}{X} + \Phi(2T) + Xgk\log T +
\Phi(T+Y_k)+g^{-k}\log T .
\end{align}

Now we take $X=1/\sqrt{\log T}$, $g=e$ and $k=\lfloor\log\log T
/2\rfloor+1$. Note that in this case we have $1\le Y_k\le e$. From
\eqref{eq bound Delta 2 X} it follows that
$$
\arg G(1/2+iT)\ll \sqrt{\log T}\log\log T+\Phi(2T),
$$
and this gives \eqref{N1 thm}.   \qed

\section*{acknowledgement}

I would like to express my gratitude to my advisor, Professor Steve
Gonek, for suggesting this problem, and for a number of inspiring
conversations, which led to the proof presented in the paper. I also
thank him for very helpful suggestions on the content and the
presentation of this paper. Thanks also go to Siegfred Baluyot for
his helpful suggestion on an earlier draft, and to Professor Yoonbok
Lee and the anonymous referee for their valuable comments.


\begin{thebibliography}{99}

\bibitem{Aka} H. Akatsuka, \textit{Conditional estimates for error terms related to the distribution of zeros of $\zeta'(s)$}, J. Number Theory 132 (2012), no. 10, 2242-2257.

\bibitem{Ber} B. C. Berndt, \textit{The number of zeros for $\zeta^{(k)}(s)$}, J. London Math. Soc.(2) 2 (1970), 577-580.

\bibitem{FGH} D. W. Farmer, S. M. Gonek and C. P. Hughes, \textit{The maximum size of L-functions}, J. Reine Angew. Math. 609 (2007), 215-236.

\bibitem{Lit} J. E. Littlewood, \textit{On the zeros of the Riemann zeta-function}, Proc. Camb. Philos. Soc. 22 (1924), 295-318.


\bibitem{Spe} A. Speiser, \textit{Geometrisches zur Riemannschen Zetafunktion}, Math. Ann. 110 (1934), 514-521.

\bibitem{Tit} E. C. Titchmarsh, \textit{The theory of the Riemann zeta-function}, 2nd ed., (ed. D. R. Heath-Brown; Oxford Science Publications, Oxford, 1986).

\bibitem{Zha} Y. Zhang, \textit{On the zeros of $\zeta'(s)$ near the critical line}, Duke Math. J. 110 (2001), 555-572.

\end{thebibliography}
\end{document}